\documentclass[12pt]{article}

\usepackage{amssymb}

\usepackage{
  latexsym,
  amsmath,
  amstext,
  amsbsy,
  amsopn,
  amsfonts
  }

\newcommand{\p}{\partial}

\newcommand{\dd}{\hspace{1pt}{\rm d}\hspace{0.5pt}}
\newcommand{\ee}{{\rm e}\hspace{1pt}}
\newcommand{\ii}{{\rm i}\hspace{1pt}}

\newcommand{\R}{\mathbb R}
\newcommand{\N}{\mathbb N}
\newcommand{\C}{\mathbb C}

\newcommand{\K}{\mathbb K}
\newcommand{\rH}{{\rm H}}

\newcommand{\cD}{{\mathcal D}}
\newcommand{\cDs}{{\mathcal D^\ast}}
\newcommand{\cE}{{\mathcal E}}
\newcommand{\cEs}{{\mathcal E^\ast}}
\newcommand{\cF}{{\mathcal F}}

\newcommand{\cT}{{\mathcal T}}

\newcommand{\Es}{{E^\ast}}

\newcommand{\qed}{\hfill $\Box$}

\newtheorem{lem}{Lemma}
\newtheorem{exmp}[lem]{Example}
\newtheorem{rem}[lem]{Remark}
\newtheorem{prop}[lem]{Proposition}
\newtheorem{cor}[lem]{Corollary}
\newtheorem{defn}[lem]{Definition}

\begin{document}

\title{On the algebraic dual of $\cD(\Omega)$}
\vspace{5mm}
\author{Michael Oberguggenberger\\
{\normalsize Institut f\"ur Grundlagen der Technischen Wissenschaften}\\
{\normalsize Arbeitsbereich f\"ur Technische Mathematik}\\
{\normalsize Universit\"at Innsbruck, A-6020 Innsbruck, Austria}\\
{\normalsize E-mail: Michael.Oberguggenberger@uibk.ac.at}}
\date{}
\maketitle

{\bf Abstract}

This paper is concerned with the algebraic dual $\cDs(\Omega)$ of the space of test
functions $\cD(\Omega)$. The emphasis is on failures and successes of $\cDs(\Omega)$
as compared to the continuous dual $\cD'(\Omega)$, the space of distributions.
Topological properties, operations with elements of $\cDs(\Omega)$ and applications
to linear partial differential equations are discussed.
\vspace{5mm}
\section{Introduction}
\label{Sectintro}
%
%
Topology plays a prominent and indispensable role in the theory of distributions,
as has been emphasized e.g. by John Horv\'ath in his monograph \cite{Hor66} as well as at numerous other places
\cite{Hor70,Hor72a,Hor72b,Hor74,Hor77a,Hor77b,Hor91}.
For example, the space of distributions $\cD'(\Omega)$ on an open subset $\Omega$ of $\R^n$
is defined as the continuous dual of the space of test functions $\cD(\Omega)$;
similarly, all other spaces of distributions can be viewed as continuous duals.
The fact that the elements of $\cD'(\Omega)$
are {\em continuous} linear functionals is essential in many constructions
as well as applications to partial differential equations. When teaching
distribution theory one usually has to spend some effort on explaining the topology
of $\cD(\Omega)$ as a locally convex inductive limit of Fr\'echet spaces. Thus
occasionally the question arises what would happen if one dropped continuity from
the definition of $\cD'(\Omega)$ and considered the algebraic dual $\cDs(\Omega)$
instead. In this paper I wish to pursue this question, and in particular, show
what goes wrong with $\cDs(\Omega)$ at the hand of a number of examples.

While these examples will clearly exhibit the lack of certain desirable properties of
$\cDs(\Omega)$ for the purpose of analysis, it is curious to note that as a
topological vector space, $\cDs(\Omega)$ has better properties than $\cD'(\Omega)$.
Not surprisingly, certain partial differential equations that do not have solutions
in $\cD'(\Omega)$ turn out to be solvable in $\cDs(\Omega)$. For example,
constant coefficient partial differential operators have solutions in $\cDs(\Omega)$
on {\em every} open subset of $\R^n$ with arbitrary members of $\cDs(\Omega)$
on the right hand side. A similar solvability result in $\cDs(\Omega)$ will be seen to hold, e.~g.,
for the Lewy equation.
This, however, is counterbalanced by the fact that one cannot say much about the
behavior of these solutions, having lost control over their analytical properties
due to arguments involving algebraic bases.

The plan of the paper is as follows. In Section \ref{Sectnotation} some basic notions
needed in the sequel are recalled.
In Section \ref{Secttopo}
properties of $\cDs(\Omega)$ as a topological vector space are collected.
Although these results are known it seemed appropriate to arrange them in the
context of the theme of the paper.
In Section \ref{Sectgenfun}, I present a number of assertions and examples
demonstrating failures (and successes) of $\cDs(\Omega)$. To my knowledge,
these considerations have not appeared in print so far. On the positive side, we will
see that derivation, multiplication by smooth functions and sheaf theoretic arguments
work well in $\cDs(\Omega)$. On the negative side, we will encounter the failure of
convolution to regularize, difficulties with the definition of tensor products and convolution,
and the lack of the notion of local order in $\cDs(\Omega)$.
Finally, Section \ref{SectPDEs}
contains some observations on solvability of partial differential equations. We dwell a bit
on the role of $P$-convexity, hypoellipticity and fundamental solutions in $\cDs(\Omega)$ there
(part of the latter results are based on joint work with T. Todorov \cite{ObeTod94}).
%
%
\section{Notation}
\label{Sectnotation}
%
%
Throughout the paper, I follow the notation of \cite{Hor66}.
The term {\em locally convex space} will refer to a locally convex Hausdorff
topological vector space over the field $\K$ ($\K$ will be either $\R$ or $\C$ in the
sequel). If the vector spaces $F, G$ form a dual system $(F,G)$
\cite[Def. 3.2.1]{Hor66}, the \mbox{\em weak-,}\ {\em Mackey-} and {\em strong topologies} on $F$
are denoted by $\sigma(F,G), \tau(F,G)$ and $\beta(F,G)$, respectively. These are
the topologies of uniform convergence on the finite subsets of $G$, on the absolutely convex,
$\sigma(G,F)$-compact subsets of $G$, and on the $\sigma(G,F)$-bounded subsets of $G$,
respectively. The {\em algebraic dual} of a vector space $E$ is the set of all linear
maps from $E$ into $\K$ and will be denoted by $\Es$. If $E$ is a locally convex
space with topology $\cT$, the {\em continuous dual} or simply {\em dual} is
the set of linear forms continuous with respect to the topology $\cT$ and will be
denoted by $E'$. It is known that $E'$ is the dual of $E$ with respect to every
locally convex topology finer than $\sigma(E,E')$ and coarser than $\tau(E,E')$
\cite[Prop. 3.5.4]{Hor66}. A locally convex space $E$ is {\em complete} if every Cauchy filter
on $E$ converges. An absolutely convex,
absorbing and closed subset of $F$ is called a barrel. The locally convex
space $E$ is called {\em barrelled}, if every barrel is a neighborhood of zero.
The family of {\em all} absolutely convex, absorbing subsets of a vector space $E$
generates the finest (i.~e., largest) locally convex topology on $E$ \cite[Ex. 2.4.3]{Hor66},
which we denote by $\cT_\ell$.

Let $\Omega$ be an open subset of $\R^n$. Then $\cE(\Omega)$ is the space of infinitely
differentiable functions on $\Omega$ with values in $\K = \C$. Equipped with the topology of
uniform convergence on compact subsets of $\Omega$, it is a complete and metrizable
locally convex space (a {\em Fr\'echet} space). The support of a smooth function
is the closure (in $\Omega$) of the set of points on which it does not vanish.
Given a compact subset $K \subset \Omega$, $\cD_K(\Omega)$ denotes the subspace of
$\cE(\Omega)$ of smooth functions with support in $K$. The union of all $\cD_K(\Omega)$
as $K$ runs through the compact subsets of $\Omega$ is the space $\cD(\Omega)$
of {\em compactly supported} smooth functions. Its genuine topology $\cT_\cD$ is the
final locally convex topology with respect to all injections $\cD_K(\Omega) \to
\cD(\Omega)$, with which it is a strict inductive limit of Fr\'echet spaces
\cite[Sect. 2.12]{Hor66}. The space of distributions on $\Omega$, $\cD'(\Omega)$,
is the continuous dual of $\cD(\Omega)$. Given $S \in \cD'(\Omega)$ and
$\varphi \in \cD(\Omega)$, the action of $S$ on $\varphi$ is denoted by
$\langle S, \varphi\rangle$. If $U$ is an open subset of $\Omega$, there is
a natural injection of $\cD(U)$ into $\cD(\Omega)$; its transpose defines
the restriction map of $\cD'(\Omega)$ to $\cD'(U)$. The support of a
distribution $S \in \cD'(\Omega)$ is the complement of the largest open set
$U$ such that the restriction of $S$ to $U$ vanishes.

The space $\cD(\Omega)$ is densely imbedded in $\cE(\Omega)$; hence the transpose
of the imbedding is injective - this way $\cE'(\Omega)$ can be viewed as a
subspace of $\cD'(\Omega)$ and in fact be identified with the space of distributions
with compact support \cite[Prop. 4.2.3]{Hor66}. Recall also that any locally integrable
function $f$ can be viewed as a distribution, given by the action
$\varphi\to \int f(x)\varphi(x)\dd x$ for $\varphi\in\cD(\Omega)$.
In particular, the space of smooth functions $\cE(\Omega)$
is contained in $\cD'(\Omega)$.

Let $K$ be a compact subset of $\Omega, L > 0, \sigma > 1$. The space
$\cD_\sigma(\Omega,K,L)$ is defined as the subspace of $\cD_K(\Omega)$ of
functions whose $p$-th partial derivatives are bounded, uniformly on $\Omega$, by a constant
times $L^{|p|}(|p|!)^\sigma$. The inductive limit of the spaces $\cD_\sigma(\Omega,K,L)$
as $K$ runs through all compact subsets of $\Omega$ and $L\to \infty$ is the
Gevrey class of order $\sigma$, $\cD_\sigma(\Omega)$. Its continuous dual
is the space of Gevrey ultradistributions of order $\sigma$, $\cD_\sigma'(\Omega)$,
see e.~g. \cite[Chap. 7, Def. 2.1]{LioMag70}.
%
%
\section{Topological properties of $\cDs(\Omega)$}
\label{Secttopo}
%
%
{\bf General properties of algebraic duals.} All results in this section are known,
but will be useful and relevant for a proper understanding of $\cDs(\Omega)$.
We begin by collecting some properties
that hold for algebraic duals in general.
Thus let $E$ be a locally convex space and let $(e_\lambda)_{\lambda \in \Lambda}$
be an algebraic basis of $E$. Then $E$ is algebraically isomorphic
with the direct sum of $|\Lambda|$ copies of $\K$ and $\Es$ with the corresponding
direct product:
\begin{equation}\label{isomorphisms}
   E \approx \K^{(\Lambda)},\quad \Es \approx \K^\Lambda.
\end{equation}
The space $\K^{(\Lambda)}$ is equipped with the finest locally convex topology making all
injections $\K^I \to \K^{(\Lambda)}$, $I$ finite, continuous. It is clear that
this topology coincides with the finest locally convex topology on $\K^{(\Lambda)}$.
Further, the product topology on $\K^{\Lambda}$ coincides
with the weak topology $\sigma(\K^{\Lambda}, \K^{(\Lambda)})$ \cite[Prop. 3.14.3]{Hor66}.
Clearly, the dual of $E$ with respect to the finest locally convex topology $\cT_\ell$
is $\Es$. Thus, if we put the finest locally convex topology $\cT_\ell$ on $E$ and the weak topology
$\sigma(\Es,E)$ on $\Es$, the isomorphisms in (\ref{isomorphisms}) are topological.
\begin{lem} \label{weakstar}
Let $E$ be a vector space. Then:\\
(a) Every $\sigma(E,\Es)$-bounded subset of $E$ is finite dimensional.\\
(b) Every subspace of $E$ is closed with respect to the topology $\sigma(E,\Es)$.
\end{lem}
{\em Proof:} (a) If $(x_n)_{n\in\N}$ is an infinite sequence of
linearly independent members of $E$, one can find
an element $x^\ast \in \Es$ such that $\langle x_n,x^\ast\rangle = n$; thus
the set $(x_n)_{n\in\N}$ is unbounded. (b) If $L$ is a subspace of $E$ and
$x \not\in L$, one can find a linear form which vanishes on $L$ and has value 1 on $x$, say.
\qed

The topology $\sigma(\Es,E)$ is the topology of uniform convergence on the
$\sigma(E,\Es)$-bounded, finite-dimensional subsets of $E$ \cite[Ex. 3.4.1]{Hor66}, while the
topology $\beta(\Es,E)$ is the topology of uniform convergence on the
$\sigma(E,\Es)$-bounded subsets of $E$. The Mackey topology $\tau(\Es,E)$ is the topology of
uniform convergence on the absolutely convex,  $\sigma(E,\Es)$-compact subsets of $E$.
Since these are $\sigma(E,\Es)$-bounded as well, Lemma \ref{weakstar}\,(a) implies
that the \mbox{weak-,}\ Mackey- and strong topology coincide on $\Es$:
\[
   \sigma(\Es,E) = \tau(\Es,E) = \beta(\Es,E).
\]
As noted above, $\Es$ is the dual of $E$ with respect to the finest locally
convex topology $\cT_\ell$. It follows from \cite[Prop. 3.5.4]{Hor66} that
\[
 \cT_\ell = \tau(E,\Es).
\]
\begin{prop} \label{complete}
Let $E$ be a vector space. Then:\\
(a) $\Es$ is complete with respect to $\sigma(\Es,E)$.\\
(b) $E$ is complete with respect to $\tau(E,\Es)$.
\end{prop}
{\em Proof:} This follows from the isomorphisms (\ref{isomorphisms}) above and
\cite[Rem. 2.11.1]{Hor66} \cite[\S 18.5.(3)]{Koe66}. \qed

{\bf Properties related to the Banach-Steinhaus theorem.} One of the important theorems
of analysis is the Banach-Steinhaus theorem. In
one of its forms, it relates equicontinuity and pointwise boundedness of continuous
linear maps. Thus let $F,G$ be locally convex spaces and consider the statement

(S1)\qquad
\parbox{10cm}{Every pointwise bounded family of continuous linear maps from
$F$ into $G$ is equicontinuous.}

The question about the maximal class of locally convex spaces $F$ such that (S1) holds for
all locally convex spaces $G$ is answered by the {\em Banach-Steinhaus theorem}; it is the
class of barrelled locally convex spaces: {\em A locally
convex space $F$ satisfies (S1) for every locally convex space $G$ if and only if
it satisfies (S1) for $G = \K$, if and only if
it is barrelled}, see e.~g. \cite[Prop. 3.6.2]{Hor66} or \cite[Thm. 3.2.3]{Obe82}.
For applications, the following corollary (see \cite[Cor. to Prop. 3.6.5]{Hor66})
of the Banach-Steinhaus theorem is important:
\begin{prop} \label{BanachSteinhaus}
Let $F$ be a barrelled locally convex space, $(x_n^\ast)_{n\in\N}$ a sequence
of continuous linear forms on $F$ which converges pointwise, that is,
$\langle x^\ast_n, x\rangle$ converges to a limit $\langle x^\ast,x\rangle$
for every $x\in F$. Then $x^\ast$ defines a continuous linear form on $F$. \qed
\end{prop}
\begin{cor} \label{barrelledness}
Let $E$ be a vector space. Then:\\
(a) $\Es$ is barrelled with respect to the topology $\sigma(\Es,E)$.\\
(b) $E$ is barrelled with respect to the finest locally convex topology $\tau(E,\Es)$.
\end{cor}
{\em Proof:} (a) Let $X$ be a family of pointwise bounded $\sigma(\Es,E)$-continuous linear
maps from $\Es$ into $\K$. By \cite[Prop. 3.2.2]{Hor66}, $X$ is a subset of $E$.
By Lemma \ref{weakstar}, $X$ is finite dimensional. Being bounded, it is also
contained in the convex hull
of finitely many points, hence equicontinuous. Thus property (S1) holds for $F=\Es$
and $G=\K$, so $\Es$ is barrelled.
(b) In the finest locally convex topology, every barrel is a neighborhood. \qed
\begin{cor}\label{sequentiallycomplete}
If $E$ is a barrelled locally convex space, then $E'$ is sequentially complete
with respect to $\sigma(E',E)$.
\end{cor}
{\em Proof:} This is an immediate consequence of Proposition \ref{BanachSteinhaus}. \qed

We now summarize these observations in the situation of the algebraic dual
of the space of test functions and
combine them with classical facts about distribution spaces.
Let $\Omega$ be an open subset of $\R^n$. For simplicity of notation, we drop
reference to the open set $\Omega$ in the expressions for the polar topologies.
Thus $\tau(\cD,\cD')$ will mean $\tau(\cD(\Omega),\cD'(\Omega))$ in the remainder
of this section, and similarly for the other spaces and polar topologies.

Clearly, the inductive limit topology $\cT_\cD$ on $\cD(\Omega)$ coincides with
the Mackey topology $\tau(\cD,\cD')$, which in turn is strictly coarser than the
topology $\cT_\ell = \tau(\cD,\cDs)$.

\begin{prop} (a) The space of test functions $\cD(\Omega)$ is complete and
barrelled with respect to both topologies $\tau(\cD,\cD')$ and $\tau(\cD,\cDs)$.\\
(b) The continuous dual $\cD'(\Omega)$ of $\cD(\Omega)$ (with respect to $\cT_\cD$) is complete
and barrelled with the topology $\beta(\cD',\cD)$ and sequentially complete
with the topology $\sigma(\cD',\cD)$.\\
(c) The algebraic dual $\cDs(\Omega)$ of $\cD(\Omega)$ is complete
and barrelled with the topology $\sigma(\cDs,\cD) = \tau(\cDs,\cD) = \beta(\cDs,\cD)$.\\
(d) $\cE(\Omega)$ is dense, but not sequentially dense in $\cDs(\Omega)$ with respect to
the topology $\sigma(\cDs,\cD)$.\\
(e) $\cD'(\Omega)$ is the completion of $\cE(\Omega)$ with regard to $\beta(\cD',\cD)$;
$\cDs(\Omega)$ is the completion of $\cE(\Omega)$ with regard to $\sigma(\cD',\cD)$.
\end{prop}

{\em Proof:} (a) The statements about $\tau(\cD,\cDs)$ follow from
Proposition \ref{complete} and Corollary \ref{barrelledness},
the statements about $\tau(\cD,\cD')$ from Cor. to Thm. 2.12.3 and
Prop. 3.6.4 in \cite{Hor66}. (b) The statements about $\beta(\cD',\cD)$ follow, for
example, by combining the assertions from Ex. 3.7.2, Prop. 3.7.6, Ex. 3.9.6,
Prop. 3.9.9 and the sentence after Cor. to Prop. 3.9.1 in \cite{Hor66}.
The sequential completeness of $\cD'(\Omega)$ with $\sigma(\cD',\cD)$ follows from
(a) and Corollary \ref{sequentiallycomplete}.
(c) This is asserted by Proposition \ref{complete} and Corollary \ref{barrelledness}.
(d) The subspace of $\cD(\Omega)$ orthogonal
to $\cE(\Omega)$ is $\{0\}$, thus $\cE(\Omega)$ is dense in $\cDs(\Omega)$ \cite[Prop. 3.3.3]{Hor66}.
The fact that $\cE(\Omega)$ is not sequentially dense in $\cDs(\Omega)$ follows from (b), the
sequential completeness of $\cD'(\Omega)$.
(e) By (d), $\cE(\Omega)$ is dense in $\cDs(\Omega)$ whose weak topology $\sigma(\cDs,\cD)$
is complete, whence the second statement. But $\cE(\Omega)$ is weakly dense in $\cD'(\Omega)$ all the
more, hence also dense in $\cD'(\Omega)$ with respect to the topology
$\beta(\cD',\cD) = \tau(\cD',\cD)$ \cite[Prop. 3.4.3]{Hor66}. \qed

{\bf Properties related to the closed graph theorem.} A second important theorem of
analysis is the closed graph theorem to which we
now turn. Thus let $F,G$ be locally convex spaces and consider the statement

(S2)\qquad
\parbox{10cm}{Every linear map from $F$ into $G$ whose graph is a closed subset
of $F\times G$ is continuous.}

The classical {\em closed graph theorem} of Banach \cite{Ban32} says that
statement (S2) is true if both $F$ and $G$ are Fr\'echet spaces. In order to
extend this theorem to more general classes of spaces, we might first fix
the class on the left hand side, say to the class of barrelled spaces.
What spaces then are admitted on the right hand side to make (S2) true?
Consider a locally convex space $G$. A subset $L \subset G'$ is called
$\nu(G',G)$-closed, if its intersections $L\cap U$ with all equicontinuous subsets $U$
of $G'$ are $\sigma(G',G)$-closed (in $U$). The locally convex space $G$ is
called {\em fully complete} or a {\em Pt\'ak space} if all $\nu(G',G)$-closed
subspaces of $G'$ are $\sigma(G',G)$-closed. The closed graph theorem of
Robertson and Robertson \cite{Rob56} states that (S2) is true if $F$ is
barrelled and $G$ a Pt\'ak space, see also \cite[Thm. 3.17.4]{Hor66}.

Let us aim directly at the class of all locally convex spaces $G$, such that
(S2) is true for every barrelled space $F$. A space $G$ with this
property is called {\em (barrelled)-minimal} or an {\em infra-(s)-space}.
The question whether
the spaces $\cD(\Omega)$ and $\cD'(\Omega)$ are Pt\'ak spaces (or more
generally infra-Pt\'ak spaces \cite[Sect. 3.10]{Hor66}), was settled
in the negative by Valdivia \cite{Val74,Val77}. Since every barrelled infra-(s)-space is
an infra-Pt\'ak space (\cite[\S34.9.(8)]{Koe79}, \cite[S. 7.3.9]{Obe82}),
it follows that $\cD(\Omega)$ and $\cD'(\Omega)$ are not (barrelled)-minimal;
neither with their strong nor their weak topologies \cite[Bsp. 7.3.10, Bsp. 7.3.11]{Obe82}.

The situation is different with the algebraic dual $\cDs(\Omega)$. Observe
first that every subspace of $\cD(\Omega)$ is $\sigma(\cD,\cDs)$-closed
(Lemma \ref{weakstar}).
It follows that $\cDs(\Omega)$ is a Pt\'ak space.
But much more is actually true: $\cDs(\Omega)$ is (locally convex)-minimal, that is,
statement (S2) holds with every locally convex space $F$ when the target
space $G$ is $\cDs(\Omega)$. In fact, a locally convex space is (locally convex)-minimal
if and only if it is isomorphic to $\K^\Lambda$ for some index set $\Lambda$.
This is a consequence of the closed graph theorem of K\={o}mura \cite{Kom62},
see e.~g. \cite[Thm. 7.2.4]{Obe82} for further details.

On the other hand, the space $\cD(\Omega)$ with the finest locally convex topology
is clearly not (barrelled)-minimal, because $\tau(\cD,\cD')$ is a strictly coarser
barrelled topology. We may summarize what has just been deduced as follows:

\begin{prop} (a) $\cDs(\Omega)$ is a Pt\'ak space and (locally convex)-minimal
in the topology $\sigma(\cDs,\cD) = \tau(\cDs,\cD) = \beta(\cDs,\cD)$.\\
(b) $\cD'(\Omega)$ is not (barrelled)-minimal, neither with respect to $\sigma(\cD',\cD)$
nor with respect to $\beta(\cD',\cD)$.\\
(c) $\cD(\Omega)$ is not (barrelled)-minimal, neither with respect to $\tau(\cD,\cD')$
nor with respect to $\tau(\cD,\cD^\ast)$. \qed
\end{prop}

With regard to statement (S2), the situation is much better for the
Fr\'echet space $\cE(\Omega)$ and its continuous dual. For example, $\cE'(\Omega)$
is a Pt\'ak space in the topology $\beta(\cE',\cE)$ \cite[Prop. 3.17.6]{Hor66}. Also, both
$\cD(\Omega)$ and $\cD'(\Omega)$ are ultrabornological (i.~e., inductive
limits of an arbitrary family of Fr\'echet spaces) and (ultrabornological)-minimal,
thanks to De Wilde's theory \cite{DeW69}. For more details on the closed
graph theorem and its historical aspects we refer to
\cite{Hor73a,Hor73b,Koe79,Obe82}.
%
%
\section{$\cDs(\Omega)$ as a space of generalized functions}
\label{Sectgenfun}
%
%
{\bf Derivation and multiplication.} Let $\Omega$ be an open subset of $\R^n$.
We begin by collecting some positive
results about the algebraic dual $\cDs(\Omega)$, to show that $\cDs(\Omega)$ may serve as a
space of generalized functions on $\Omega$.
First, an element $S\in\cDs(\Omega)$ can be differentiated and multiplied by smooth
functions. The definitions follow the same lines as in \cite[Sect. 4.3, Sect. 4.6]{Hor66}.
Thus let $p \in \N_0^n$. Then the $p$-th partial derivative of $S$ is defined
as
\[
   \langle \p^p S,\varphi\rangle = (-1)^{|p|}\langle S, \p^p \varphi\rangle
\]
and the product of $S$ with a smooth function $\alpha\in \cE(\Omega)$ by
\[
   \langle \alpha S,\varphi\rangle = \langle S, \alpha\varphi\rangle
\]
for $\varphi \in \cD(\Omega)$. For the one-dimensional case $\Omega = \R$,
a repetition of the classical proof shows:
\begin{prop}
The map $\p : \cDs(\R) \to \cDs(\R)$ is surjective; its kernel consists of the
one-dimensional subspace of constant functions.
\end{prop}
{\em Proof:} Following \cite[Sect. 4.3]{Hor66}, we denote the image of
$\p : \cD(\R) \to \cD(\R)$ by $H$. A test function $\chi$ belongs to $H$
if and only if its integral vanishes. Taking an arbitrary test function $\varphi_0$
with integral 1, every $\varphi \in \cD(\R)$ can be written as
\[
   \varphi = \lambda\varphi_0 + \chi
\]
with $\chi = \p\psi\in H$ and $\lambda = \int_{-\infty}^{\infty} \varphi(x)\dd x$.
Let $T\in \cDs(\R)$ and $\p T = 0$. Then
\[
   \langle T,\varphi\rangle = \lambda\langle T,\varphi_0\rangle - \langle \p T,\psi\rangle
     = \langle T,\varphi_0\rangle\int_{-\infty}^{\infty} \varphi(x)\dd x.
\]
This shows that the action of $T$ is given by the constant function $\langle T,\varphi_0\rangle$.
On the other hand, we have the direct sum decomposition
\[
   \cD(\R) = \C\,\varphi_0 \oplus H,
\]
and $-\p :\cD(\R) \to H$ is injective. Thus given $S \in \cDs(\R)$,
the element $T \in \cDs(\R)$,
\[
   \langle T,\varphi\rangle = -\langle S,\psi\rangle
\]
is well defined (where $\psi$ has the same meaning as above), and clearly
$\p T = S$. \qed

Let $h\in\R^n$. The translate of a function $f\in \cE(\R^n)$ is defined by
$(\tau_h f)(x) = f(x-h)$. The translate of an element $T$ of $\cDs(\R^n)$
can be defined by
\[
   \langle \tau_hT,\varphi\rangle = \langle T, \tau_{-h}\varphi\rangle
\]
as in the case of distributions \cite[Def. 4.3.2]{Hor66}.
However, the map
$h \to \tau_h T$ from $\R^n$ to $\cDs(\R^n)$ is not continuous
when $\cDs(\Omega)$ is equipped with the topology $\sigma(\cDs,\cD)$. This will follow from
Example\;\ref{exampleconvolution} below.
Nevertheless, we may state as in \cite[Def. 4.3.3]{Hor66} that an element $T \in \cDs(\R^n)$
is {\em independent of the variable} $x_j$ if $\tau_h T = T$ for all
vectors $h\in \R^n$ parallel to the $x_j$-axis.

{\bf Supports and restrictions.} We now turn to the sheaf properties of $\cD^\ast(\Omega)$. First, if $U$ is an open
subset of $\Omega$ every function $\varphi \in \cD(U)$ can be considered as
an element of $\cD(\Omega)$, extending it by 0 outside $U$. Thus an element
$T \in \cDs(\Omega)$ can be restricted to $U$ by the prescription
\[
   \langle T|U, \varphi\rangle = \langle T, \varphi\rangle
\]
for $\varphi \in \cD(U)$. Clearly, if $V\subset U$ then $T|V = (T|U)|V$.
We say that $T$ vanishes on $U$ if $T|U = 0$.
Given an open cover $(\Omega_\iota)_{\iota\in I}$ of $\Omega$ we have
that $T$ vanishes on $\Omega$ if and only if all its restrictions to $\Omega_\iota$
vanish for every $\iota\in I$. Indeed, there is a locally finite, infinitely
differentiable partition of unity $(\alpha_\iota)_{\iota\in I}$
subordinated to the cover $(\Omega_\iota)_{\iota\in I}$ \cite[Thm. 2.12.4]{Hor66}.
For $\varphi\in \cD(\Omega)$ we thus have
\[
   \langle T,\varphi\rangle = \sum_{\iota\in I}\langle T,\alpha_\iota\varphi\rangle = 0
\]
since the sum contains only finitely many terms when $\varphi$ is fixed.
In the same vein, given a family of elements $T_\iota \in \cDs(\Omega_\iota),
\iota \in I$, such that $T_\iota = T_\kappa$ on each non-empty common domain
$\Omega_\iota \cap \Omega_\kappa$, there is a unique element $T\in \cDs(\Omega)$
such that $T|\Omega_\iota = T_\iota$ for every $\iota\in I$. The proof of this
fact is the same as in \cite[Prop. 4.2.4]{Hor66}; actually shorter since the
continuity argument is not needed. We have proven:
\begin{prop} The assignment $\Omega \to \cDs(\Omega)$ defines a sheaf of
locally convex spaces on $\R^n$. \qed
\end{prop}
In particular, the support of an element $T\in \cDs(\Omega)$ is well defined
as the complement of the largest open subset of $\Omega$ on which it vanishes.

Here comes the first major difference of the behavior of $\cDs$ as compared
to $\cD'$. As mentioned in Section \ref{Sectnotation}, the elements of $\cE'(\Omega)$
can be identified with the compactly supported distributions. This is no longer
 the case in the setting of the algebraic duals: $\cD(\Omega)$ is not a
dense subspace of $\cE(\Omega)$ with respect to the finest locally convex
topology $\tau(\cE,\cEs)$, but rather a closed subspace (Lemma \ref{weakstar});
hence the transpose of this
imbedding is not an injective map from $\cEs(\Omega)$ to $\cDs(\Omega)$.
On the contrary, we have an injection in the reverse direction. To see this,
let $N$ be an algebraic supplement of $\cD(\Omega)$ in $\cE(\Omega)$. The
map
\begin{equation}\label{injection}
   i: \cDs(\Omega) \to \cEs(\Omega),\quad \langle i(T),\varphi\rangle = \langle T,\psi\rangle
\end{equation}
where $T\in \cDs(\Omega)$ and $\varphi = \psi + \chi$ with $\psi \in \cD(\Omega)$ and
$\chi \in N$ is clearly linear and injective. This way $\cDs(\Omega)$ becomes
a subspace of $\cEs(\Omega)$, and membership in $\cEs(\Omega)$ does not correspond to
any support property. There are many elements of $\cEs(\Omega)$ with the same
action as $T$ on $\cD(\Omega)$, namely all those of the form
$\varphi \to \langle T,\psi\rangle + \langle T', \chi\rangle$ where
$T'$ is some linear functional on $N$.

A similar situation arises with respect to the spaces of distributions of finite
order. Let $\cD^m(\Omega), m \in \N_0$, denote the space of $m$-times continuously
differentiable functions with compact support. Its continuous dual ${\cD'}^m(\Omega)$
is the space of distributions of order (at most) $m$ and is a subspace of
$\cD'(\Omega)$. Again, in the setting of algebraic duals, the injections are
reversed: Letting $N^m$ be an algebraic supplement of $\cD(\Omega)$ in $\cD^m(\Omega)$,
the same reasoning as in (\ref{injection}) leads to a linear injection of
$\cDs(\Omega)$ in $\cDs^m(\Omega)$. Thus the notion of {\em order} has no
meaning for the elements of $\cDs(\Omega)$. Indeed, we will shortly exhibit
elements that do not arise as distributions of locally finite order.
\begin{exmp}\label{exampleinfiniteorder}
Let $M$ be the subspace of $\cD(\R)$ of test functions whose sequence of derivatives
at zero is summable:
\[
  M = \{\psi \in \cD(\R): \sum_{k=0}^\infty|\p^k\psi(0)| < \infty\},
\]
and let $N$ be an algebraic supplement of $M$ in $\cD(\R)$. The prescription
\[
   \langle T,\varphi\rangle = \sum_{k=0}^\infty\p^k\psi(0)
\]
where $\varphi = \psi + \chi$ with $\psi \in M, \chi \in N$ defines an
element $T \in \cDs(\R)$. Involving infinitely many derivatives at zero, $T$ is
not a continuous functional on $\cD(\R)$ with respect to $\cT_D$,
thus does not belong to $\cD'(\R)$.
\end{exmp}
Actually, the simple algebraic argument in (\ref{injection}) can be
generalized to show that the spaces of Gevrey ultradistributions are also contained
as subspaces of $\cDs(\Omega)$. Thus we have the somewhat curious sequence of
inclusions $(\sigma > 1)$
\[
 \cE'(\Omega) \subset \cD'(\Omega) \subset \cD'_\sigma(\Omega) \subset
   \cD^\ast_\sigma(\Omega) \subset \cDs(\Omega) \subset \cEs(\Omega).
\]
Example \ref{exampleinfiniteorder} is also of interest from the viewpoint of
supports: the support of the distribution $T$ defined there
is $\{0\}$. Indeed, if $\varphi \in \cD(\R\setminus\{0\})$, then
$\varphi$ belongs to $M$ and $\langle T,\varphi\rangle = 0$.
We see that an element of $\cD(\R)$ whose support is $\{0\}$ need not
be a finite linear combination of the Dirac measure and its derivatives
(as opposed to the case of distributions, \cite[Prop. 4.4.5]{Hor66}).
Here is another example of this phenomenon.
\begin{exmp}
Fix an element $\varphi_0$ of $\cD(\R)$ such that $\varphi_0(0) = 1$ and
$\p^k\varphi_0(0) = 0$ for all $k\geq 1$. Let $H$ be an algebraic supplement
of the one-dimensional space $\C\,\varphi_0$ in $\cD(\R)$. Define
a linear form $S$ on $\cD(\R)$ by
\[
   \langle S,\varphi\rangle = \lambda,
\]
where $\varphi = \lambda \varphi_0 + \chi$ with $\lambda \in \K$ and $\chi \in H$.
Then $S$ is not a finite linear combination of the Dirac measure $\delta$
and its derivatives. Indeed, assume to the contrary that $S = \sum_{p=0}^ma_p\p^p\delta$ for
some $m \in \N$ and certain coefficients $a_p$. Letting $\varphi = \lambda \varphi_0 + \chi$
with $\chi \in H$, we would have that
\begin{equation} \label{examplesupport}
\langle \sum_{p=0}^ma_p\p^p\delta,\varphi\rangle
  = \sum_{p=0}^m(-1)^pa_p\p^p\varphi(0) = \lambda + \sum_{p=0}^m(-1)^pa_p\p^p\chi(0).
\end{equation}
If this expression represented $\langle S,\varphi\rangle$ it should equal $\lambda$, for arbitrary
choices of $\chi \in H$. This is not the case, because one can always find elements
$\chi$ of $\cD(\R)$ which are not multiples of $\varphi_0$ such that the sum on the
right hand side of (\ref{examplesupport}) does not vanish.
\end{exmp}

{\bf Convolutions and tensor product.} We now arrive at a more severe failure of $\cD(\R)$, and that is the failure
of convolutions to regularize.
\begin{lem}\label{lemmaindependence}
Let $\varphi$ be a nonzero element of $\cD(\R)$. Then the family of translates
$(\tau_{h}\varphi)_{h\in\R}$ is linearly independent in $\cD(\R)$.
\end{lem}\label{linearindependent}
{\em Proof:} Assume that
\[
    \sum_{p=0}^m a_p \tau_{h_p}\varphi(x) \equiv 0 \ \mbox{on}\ \R
\]
for certain $m\in \N, h_p\in\R$ and $a_p \in \C$. Taking the Fourier transform of this
equation, we have that
\[
    \Big(\sum_{p=0}^m a_p \ee^{-\ii h_p\xi}\Big)\cF\varphi(\xi) \equiv 0 \ \mbox{on}\ \R.
\]
Since both factors can be extended as entire functions of $\xi$ to the complex
plane and the ring of holomorphic functions has no zero divisors, it
follows that
\[
   \sum_{p=0}^m a_p \ee^{-\ii h_p\zeta} \equiv 0 \ \mbox{on}\ \C.
\]
But exponentials of different phase are linearly independent, thus
all coefficients $a_p$ necessarily vanish. \qed
\begin{defn}\label{defconvolution}
Let $S \in \cDs(\R^n), \varphi\in \cD(\R^n), x\in\R^n$. The convolution
of $S$ and $\varphi$ at the point $x$ is defined as
$S\ast\varphi(x) = \langle S, \tau_x\check{\varphi}\rangle$.
\end{defn}
Here $\check{\varphi}(y) = \varphi(-y)$; in abusive notation involving
the independent variable inside the duality brackets, the definition
may become intuitively clearer:
\[
   S\ast\varphi(x) = \langle S(y), \varphi(x-y)\rangle.
\]
As in the case of distributions, the convolution with a test function yields a
function from $\R^n$ to $\K$. However, it need no longer be smooth, not even continuous.
\begin{exmp}\label{exampleconvolution}
Let $\varphi$ be a nonzero element of $\cD(\R)$.
Consider the sequence $h_n = 1/n, h_0 = 0$ in $\R$. By Lemma \ref{linearindependent}
the sequence of translates $(\tau_{h_n}\check{\varphi})_{n\in\N_0}$ is a linearly independent
subset of $\cD(\R)$. Denote by $M$ its span and by $N$ an algebraic supplement of $M$
in $\cD(\R)$. Define an element $S$ of $\cDs(\R)$ by
\[
   \langle S, \check{\varphi} \rangle = 0,\quad
   \langle S, \tau_{h_n}\check{\varphi} \rangle = n\ {\rm for}\ n\geq 1,\quad
   \langle S, \chi \rangle = 0\ {\rm for}\ \chi\in N.
\]
Then obviously
\[
  S\ast\varphi\,\big(\tfrac{1}{n}\big) = \langle S, \tau_{h_n}\check{\varphi} \rangle = n,
      \quad S\ast\varphi\,(0) = \langle S, \check{\varphi} \rangle = 0
\]
so that the function $x \to S\ast\varphi(x)$ is discontinuous at zero.
\end{exmp}
When $S$ belongs to $\cD'(\R^n)$ and $\varphi$ to $\cD(\R^n)$,
the map $x \to S\ast\varphi(x)$ is smooth. As can be seen from \cite[Prop. 4.10.1]{Hor66},
the fact that $S$ is a {\em continuous} functional on $\cD(\R^n)$ with
respect to the topology $\cT_\cD$ is at the core of the proof of this property.

Similar difficulties arise when one wants to define the tensor product in the
setting of the algebraic duals. Thus let $\Xi$ be an open subset of $\R^k$ and
$\rH$ an open subset of $\R^l$. Let $S \in \cD'(\Xi), T \in \cD'(\rH)$.
Given $\chi \in \cD(\Xi\times\rH)$, the map (notation as explained
after Definition \ref{defconvolution})
\begin{equation} \label{tensorformula}
   x \to \langle T(y), \chi(x,y) \rangle
\end{equation}
belongs to $\cD(\Xi)$ \cite[Lemma 4.8.2]{Hor66}. Thus one may define the
tensor product of $S$ and $T$ by
\[
   \langle S\otimes T, \chi \rangle = \langle S(x), \langle T(y), \chi(x,y) \rangle\rangle
\]
and it belongs to $\cD'(\Xi\times\rH)$ \cite[Lemma 4.8.3]{Hor66}. Alternatively, let us denote
by $\cD(\Xi)\otimes\cD(\rH)$ the span in $\cD(\Xi\times\rH)$ of elements of
the form $\chi(x,y) = (\varphi\otimes\psi)(x,y) = \varphi(x)\psi(y)$. Then $\cD(\Xi)\otimes\cD(\rH)$
is dense in $\cD(\Xi\times\rH)$ \cite[Prop. 4.8.1]{Hor66}, and $S\otimes T$ turns out
to be the unique distribution $R\in \cD'(\Xi\times\rH)$ such that
\[
   \langle R, \varphi\otimes\psi \rangle = \langle S, \varphi\rangle\langle T, \psi\rangle,
\]
see \cite[Def. 4.8.1]{Hor66}. Both approaches fail in the setting of the algebraic duals.
\begin{exmp}
Similar to Example \ref{exampleconvolution}, fix an element $\chi \in \cD(\R\times\R)$
such that $\chi(x,y) = 1$ when $\max(|x|,|y|)\leq 1$. The function
$(x,y) \to \ee^{\ii xy}\chi(x,y)$ belongs to $\cD(\R\times\R)$ as well, and the
family of functions
\[
  y \to \chi(0,y),\quad y \to \ee^{\ii\frac{1}{n}y}\chi(\tfrac{1}{n},y),\ n \in  \N,
\]
is a linearly independent subset of $\cD(\R)$. Indeed,
for $|y|\leq 1$ we actually have that $\ee^{\ii\frac{1}{n}y}\chi(\tfrac{1}{n},y)
\equiv \ee^{\ii\frac{1}{n}y}$, and this family is linearly independent, as noted
in the proof of Lemma \ref{lemmaindependence}. Now define $T\in\cDs(\R)$ by
\[
   \langle T(y), \chi(0,y)\rangle = 0,\quad
   \langle T(y), \ee^{\ii\frac{1}{n}y}\chi(\tfrac{1}{n},y)\rangle = n
\]
for $n\in \N$ and extend it by zero on an algebraic supplement of the span of this family of
functions. Then the map required in (\ref{tensorformula}), with $\ee^{\ii xy}\chi(x,y)$
in place of $\chi(x,y)$,
\[
   f : x \to \langle T(y), \ee^{\ii xy}\chi(x,y)\rangle
\]
is again discontinuous at $x=0: f(\frac{1}{n}) = n, f(0) = 0$.
\end{exmp}
As an alternative approach to defining the tensor product of elements
$S \in \cD'(\Xi), T \in \cD'(\rH)$ one could start
with the subspace $\cD(\Xi)\otimes\cD(\rH)$ of $\cD(\Xi\times\rH)$ and
define
\[
   \langle R, \sum_{i,j} a_{ij}\varphi_i\otimes\psi_j \rangle
   = \sum_{i,j} a_{ij}\langle S,\varphi_i\rangle\langle T,\psi_j \rangle.
\]
But -- as every subspace of $\cD(\Xi\times\rH)$ -- the subspace $\cD(\Xi)\otimes\cD(\rH)$
is closed for the finest locally convex topology. Thus there are many extensions
of $R$ to all of $\cD(\Xi\times\rH)$, and so it remains ambiguous how to define
$S\otimes T$ in this way (and to keep control about, e.~g., consistency with
classical definitions).

As presented in \cite[Sect. 4.9]{Hor66}, the definition of the convolution
of two distributions (with supports in favorable position) is based on the
tensor product of distributions. From what has just been said, it is impossible to give a meaning to
the convolution of two elements of $\cDs(\R^n)$ along these lines.
However, one could try to define the convolution of an element $S$ of
$\cDs(\R^n)$ with a distribution $T \in \cD'(\R^n)$ as follows.
Recall that the inflection of $T$ is defined by
\[
   \langle \check{T},\varphi \rangle = \langle T,\check{\varphi}\rangle.
\]
Given $\psi \in \cD(\R^n)$, the convolution of the distribution $\check{T}$ with
$\psi$ is a well defined smooth function, that is, $\check{T}\ast \psi$ belongs to
$\cE(\R^n)$ if $T\in \cD'(\R^n)$ and to $\cD(\R^n)$ if $T\in \cE'(\R^n)$
\cite[Prop. 4.10.1, Prop. 4.9.2]{Hor66}.
\begin{defn}\label{generalizedconvolution}
Let $S \in \cEs(\R^n)$ and $T\in \cD'(\R^n)$. Then the convolution $S\star T$
is defined as an element of $\cDs(\R^n)$ by
\begin{equation}   \label{genconvolutionformula}
  \langle S\star T,\psi\rangle = \langle S, \check{T}\ast\psi\rangle
\end{equation}
for $\psi \in \cD(\R^n)$.
\end{defn}
\begin{rem} \label{generalizedconvolutioncompact}
By what has been said just before Definition \ref{generalizedconvolution},
the right hand side of
formula (\ref{genconvolutionformula}) makes sense if $S\in \cDs(\R^n)$
and $T\in \cE'(\R^n)$. On the other hand, $\cDs(\R^n)$ is imbedded in $\cEs(\R^n)$ by means
of the injection $i$ given in (\ref{injection}), and so one may also consider
the convolution $i(S)\star T$ according to Definition \ref{generalizedconvolution}.
Due to the construction of the injection $i$, the two formulas give rise to
the same result when $T\in \cE'(\R^n)$, because $\langle i(S),\varphi\rangle
= \langle S,\varphi\rangle$ when $\varphi$ has compact support.
\end{rem}
However, other than that not much can be said about consistency with classically
defined convolutions. In fact, Definition \ref{generalizedconvolution} is not
consistent with Definition \ref{defconvolution} when $T$ is a test function itself.
\begin{exmp}
We continue with Example \ref{exampleconvolution}. First observe that if $x\in \R$
is not one of the members of the sequence $h_n, n\in \N_0$, then $\tau_x\check{\varphi}$
does not lie in its span $M$ (Lemma \ref{linearindependent}). In addition, if we take
$\varphi$ with support in the half line $(-\infty,0]$, then the function
$\varphi\ast\check{\varphi}$ does not belong to $M$ either (because it is symmetric,
while all members of M vanish on $(-\infty,0]$). Thus we may modify the
direct sum composition $\cD(\R) = M\oplus N$ as follows: we adjoin $\varphi\ast\check{\varphi}$
to $M$ and set up the algebraic complement $N$ in such a way that each $\tau_x\check{\varphi}$
belongs to $N$ when $x$ is not equal to one of the members $h_n$. We also modify
the definition of $S\in\cDs(\R)$ as follows:
\[
   \langle S, \check{\varphi} \rangle = 0,\quad
   \langle S, \tau_{h_n}\check{\varphi} \rangle = n\ {\rm for}\ n\geq 1,\quad
   \langle S,\varphi\ast\check{\varphi}\rangle = 1,\quad
   \langle S, \chi \rangle = 0\ {\rm for}\ \chi\in N.
\]
As a consequence, we have that $S\ast\varphi(\frac{1}{n}) = n$, while
$S\ast\varphi(x) = 0$ for all other $x\in\R$; in particular, the function
$x\to S\ast\varphi(x)$ is zero almost everywhere. If we choose to view it
as an element of $\cDs(\R)$ by means of the imbedding of $\cD'(\R)$, it
is the zero element. On the other hand, the convolution of $S$ and $\varphi$
according to Definition \ref{generalizedconvolution} is given by
$\langle S\star\varphi, \psi\rangle = \langle S, \check{\varphi}\ast \psi\rangle$
for $\psi\in \cD(\R)$. Taking in particular $\psi = \varphi$ we have by construction
$\langle S\star\varphi, \varphi\rangle = 1$, clearly inconsistent with
Definition \ref{defconvolution} according to which
$\langle S\ast\varphi, \varphi\rangle = \int S\ast\varphi(x)\,\varphi(x)\,\dd x = 0$.
\end{exmp}
Thus serious problems arise when one attempts to define convolutions in $\cDs(\R^n)$.
But it is worthwhile to note that the convolution introduced in Definition \ref{generalizedconvolution}
behaves well with respect to derivatives. If $S,T$ are as in Definition
\ref{generalizedconvolution} or in Remark \ref{generalizedconvolutioncompact}
then
\begin{equation}  \label{convolutionderivative}
   \p_j(S\star T) = (\p_j S)\star T = S\star (\p_j T).
\end{equation}
This follows immediately from formula (\ref{genconvolutionformula}), the corresponding
property of convolution of distributions, and the definition of partial derivatives
on $\cDs(\R^n)$.
%
%
\section{Solving linear partial differential equations in $\cDs(\Omega)$}
\label{SectPDEs}
%
%
We begin this section by an elementary observation on surjections of
algebraic duals.
\begin{prop} \label{propositionsurjection}
Let $E$ be a vector space and $P: \Es \to \Es$ a linear mapping.
Then $P$ is surjective if and only if its transpose $^tP:E\to E$ is injective.
\end{prop}
{\em Proof:} Assume first that $^tP$ is injective. Let $N$ be an algebraic supplement of $^tP(E)$ in $E$. Given $y^\ast
\in \Es$, define an element $x^\ast \in \Es$ by $\langle x^\ast,z\rangle =
\langle y^\ast ,y\rangle$ for $z =\ ^tP(y) \in\ ^tp(E)$, $\langle x^\ast,z\rangle =
0$ for $z \in N$. Since $^tP$ is injective, the element $x^\ast$ of $\Es$ is well
defined, and clearly $\langle P(x^\ast), x\rangle = \langle x^\ast,\, ^tP(x)\rangle =
\langle y^\ast, x\rangle$ for $x\in E$.
Conversely, assume that $P$ is surjective and let $^tP(y) = 0$. Then
$\langle y^\ast,y\rangle = \langle x^\ast, ^tP(y) \rangle = 0$ for all
$y^\ast \in \Es$, taking $x^\ast \in \Es$ such that $y^\ast = P(x^\ast)$.
Thus $y = 0$. \qed

We shall now explore what this purely algebraic result (no continuity is needed)
can or cannot say about solvability of partial differential equations. Thus
let $\Omega$ be an open subset of $\R^n$ and
\begin{equation*}
   P(x,\p) = \sum_{|p| \leq m}\, a_p (x)\, \p^p
\end{equation*}
be a linear partial differential operator with smooth coefficients
$a_p \in \cE(\Omega)$. Viewing $P(x,\p)$ as an operator acting on
$\cDs(\Omega)$, its transpose $^t P(x,\p)$ is
given by
\[ ^t P(x,\p)\,\varphi(x) =
                    \sum_{|p| \leq m}\,(-1)^{|p|}\,
                    \p^p \left( a_p (x)\varphi(x) \right)\,.
\]
By the proposition above, if $^tP(x,\p):\cD(\Omega) \to \cD(\Omega)$ is
injective, then the equation
\begin{equation}\label{PDE}
   P(x,\p)U = F
\end{equation}
has a solution $U\in \cDs(\Omega)$ for whatever $F\in \cDs(\Omega)$.
This line of arguments to establish solvability was first used by Todorov
\cite{Tod95} in the context of a factor space of the space of nonstandard
internal smooth functions. The $\cDs$-version was elaborated jointly
with him in an unpublished manuscript \cite{ObeTod94}.
\begin{exmp} \label{exampleinjective}
(a) If $\Omega$ is an arbitrary open subset of $\R^n$ and $P = P(\p)$ is a linear
partial differential operator with constant coefficients, then $^tP:\cD(\Omega) \to \cD(\Omega)$ is injective.
To see this, it suffices to take the Fourier transform of the equation $^tP(\p)\varphi = 0$ and
invoke the analyticity of $\cF\varphi$.\\
(b) If $\Omega$ is an arbitrary open subset of $\R^n$ and $P = P(x,\p)$ is a linear
partial differential operator of constant strength with analytic coefficients
then $^tP(x,\p):\cD(\Omega) \to \cD(\Omega)$ is injective (see the discussion before
Thm. 13.5.2 in \cite{Hoer83}). This is true, in particular,
when $P(x,\p)$ is an elliptic operator with analytic coefficients.\\
(c) If $\Omega$ is an arbitrary open subset of $\R^3$ and $P$
the Lewy operator
\begin{equation} \label{Lewy}
   P(x,\p) = -\frac{\p}{\p x_1}  - \ii\frac{\p}{\p y_1}
      + 2\ii(x_1 + i y_1)\,\frac{\p}{\p x_2}
\end{equation}
then $^tP(x,\p) = - P(x,\p):\cD(\Omega) \to \cD(\Omega)$ is injective.
This can be seen, e.~g., by the following elementary argument. Perform
a partial Fourier transform of the equation
$-P\varphi(x_1, y_1, x_2) \equiv 0$ with respect to the variable $x_2$.
Then $(\cF_{x_2\to z_2}\varphi)(x_1,x_2,z_2) = \psi(x_1,x_2,z_2)$ is an entire function of $z_2$ at fixed
$(x_1,y_1)$. Viewing $(x_1,y_1)$ as the complex variable $z_1 = x_1+\ii y_1$, we see that $\psi$
satisfies the equation
\begin{equation}\label{transposeLewy}
  \Big(\frac{\p}{\p \overline{z_1}} + z_1z_2\Big) \psi(z_1,z_2) \equiv 0.
\end{equation}
Setting $z_2 = 0$, (\ref{transposeLewy}) implies that the function $z_1 \to \psi(z_1,0)$ is
analytic; having compact support,
it necessarily vanishes identically. Successively differentiating (\ref{transposeLewy}) with respect to
$z_2$ and setting $z_2=0$, we observe that $\p_{z_2}^k\psi(z_1,0) = 0$ for all $z_1\in \C, k \in \N_0$.
Recalling the analyticity of $z_2\to \psi(z_1,z_2)$, it follows that $\psi$, and hence $\varphi$,
vanishes identically.\\
(d) Of course, there are many operators which are not injective on $\cD(\Omega)$. For example,
the operator $(1-x^2)^2\p_x + 2x:\cD(\R)\to\cD(\R)$ is not injective; the smooth function which equals
$\exp(-1/(1-x^2))$ for $|x|<1$ and $0$ otherwise is in its kernel. The operator
$x_2\p_{x_1} - x_1\p_{x_2}:\cD(\R^2)\to\cD(\R^2)$ is not injective; all rotationally invariant functions
belong to its kernel. There is even a fourth order elliptic operator with smooth coefficients
which is not injective from $\cD(\R^3)\to\cD(\R^3)$ \cite[Thm. 13.6.15]{Hoer83}.
\end{exmp}
We now discuss some special cases of the examples just mentioned in more detail.
The general solvability assertion in $\cDs(\Omega), \Omega\subset\R^3$, in
Example \ref{exampleinjective}\,(c) is curious in view of the fact that the operator
$P(x,\p)$ from (\ref{Lewy}) provided the first example, due to Lewy \cite{Lew57}, of an operator
with smooth coefficients which is not locally solvable in the sense of distributions.
That is, there exist smooth functions $F \in  \cE(\R^3)$ such that the equation
$P(x,\p)U = F$ does not have a solution in $\cD'(\Omega)$ for whatever open subset
$\Omega\subset\R^n$. In contrast, it does have solutions in $\cDs(\Omega)$.

Example \ref{exampleinjective}\,(a) and Proposition \ref{propositionsurjection} immediately
imply the following assertion.
\begin{cor} \label{corollarysurjective}
 Let $\Omega$ be an arbitrary open subset of $\R^n$, $P = P(\p)$ a nonzero linear
partial differential operator with constant coefficients. Then
$P(\p): \cDs(\Omega) \to \cDs(\Omega)$ is surjective. \qed
\end{cor}
This is in contrast with the classical situation that
the equation $P(\p)U = F$ has a solution $U \in \cD'(\Omega)$ for every
$F \in \cD'(\Omega)$ if and only if $\Omega$ is
{\em $P$-convex for supports as well as singular supports} \cite[Cor. 10.7.10]{Hoer83},
while the equation $P(\p)U = F$
has a solution $U \in \cD'(\Omega)$ for every $F \in \cE(\Omega)$
if and only if $\Omega$ is
{\em $P$-convex for supports} \cite[Thm. 10.6.6, Cor. 10.6.8]{Hoer83}.
In this context, we can observe that the property of {\em hypoellipticity} may be
lost when admitting solutions in $\cDs(\Omega)$. Recall that an operator
$P(\p)$ is hypoelliptic, if for whatever open subset $\Omega \subset \R^n$ and $U \in \cD'(\Omega)$,
$P(\p)U \in \cE(\Omega)$ implies $U \in \cE(\Omega)$. Let $P(\p)$ be a hypoelliptic
operator which is not elliptic. Then there exists an open subset $\Omega$ of $\R^n$ which
is not $P$-convex for supports \cite[Cor. 10.8.2]{Hoer83}. Thus there is $F \in \cE(\Omega)$ such
that the equation $P(\p)U = F$ has no solution $U\in \cD'(\Omega)$. However, by the corollary
above, it does have a solution $U\in \cDs(\Omega)$. By what has just been said, this solution
does not belong to $\cE(\Omega)$; thus solutions in $\cDs(\Omega)$ of hypoelliptic operators with
smooth right hand side need not be smooth.

A {\em fundamental solution} of a constant coefficient partial differential operator
$P(\p)$ is an element $S$ of $\cDs(\R^n)$ such that $P(\p)S = \delta$, the Dirac measure.
Corollary \ref{corollarysurjective} implies, in particular, that every nonzero
constant coefficient partial differential operator possesses a fundamental solution
in $\cDs(\R^n)$. Due to the theorem of Malgrange and Ehrenpreis
\cite{Ehr56,Mal55}, every nonzero constant coefficient partial differential operator
actually has a fundamental solution in $\cD'(\R^n)$. We refer to
\cite{OrtWag96,OrtWag97,Wagner2009} for an elegant explicit construction and a historical survey.
Proposition \ref{propositionsurjection} is just the simple portion of the Malgrange-Ehrenpreis
theorem, the difficult part of course being to prove the continuity of the
functional defined on the range of $^tP$.
\begin{prop}
(a) Let $S$ be a fundamental solution of $P(\p)$ in $\cDs(\R^n)$ and let $F\in \cE'(\R^n)$.
Then the element $U = S\star F \in \cDs(\R^n)$ is a solution of the equation
$P(\p)U = F$ in $\cDs(\R^n)$.\\
(b) Let $T$ be a fundamental solution of $P(\p)$ in $\cD'(\R^n)$ and let $G\in \cEs(\R^n)$.
Then the element $V = G\star T \in \cDs(\R^n)$ is a solution of the equation
$P(\p)V = G$ in $\cDs(\R^n)$.
\end{prop}
{\em Proof:} By Remark \ref{generalizedconvolutioncompact}, both $U$ and $V$ are well
defined elements of $\cDs(\R^n)$ according to
formula (\ref{genconvolutionformula}). Using (\ref{convolutionderivative}) we have
in case (a) that
\[
   P(\p)U = (P(\p)S)\star F = \delta\star F = F.
\]
The latter equality follows from
\[
  \langle\delta\star F,\psi\rangle
    = \langle\delta, \check{F}\ast\psi\rangle
    = (\check{F}\ast\psi)(0) = \langle F,\psi\rangle.
\]
for $\psi\in \cD(\R^n)$. In case (b), we have that
\[
   P(\p)V = G\star (P(\p)T) = G\star\delta = G,
\]
using that $\langle G\star \delta,\psi\rangle = \langle G, \check{\delta}\ast\psi\rangle
= \langle G,\psi\rangle$. \qed

Thus the simple tool of Proposition \ref{propositionsurjection} allows
to solve constant coefficient partial differential equations in $\cDs(\R^n)$,
though not much can be inferred about the properties of these solutions in general.

\end{document}